\documentclass{article}
\usepackage[margin = 1.25in]{geometry}

\usepackage{amsthm,amsmath,amsfonts,amssymb}
\usepackage{graphicx}
\usepackage{mathtools}
\usepackage{multirow}
\usepackage{epsfig,url}
\usepackage[utf8]{inputenc}

\usepackage{subfiles}

\usepackage{natbib}

\numberwithin{equation}{section}
\theoremstyle{plain}
\newtheorem{theorem}{Theorem}[section]
\newtheorem{lemma}[theorem]{Lemma}
\newtheorem{corollary}[theorem]{Corollary}
\theoremstyle{remark}
\newtheorem{remark}[theorem]{Remark}
\newtheorem{definition}[theorem]{Definition}
\newcommand{\paperTitle}{An empirical process framework for covariate balance in causal inference}

\renewenvironment{abstract}%
{%
  \vskip 0.075in%
  \centerline%
  {\large\bf Abstract}%
  \vspace{0.5ex}%
  \begin{quote}%
}
{
  \par%
  \end{quote}%
  \vskip 1ex%
}

\begin{document}

\title{\huge \paperTitle}
\date{}
\author{
\textbf{Efrén Cruz Cortés} \\
\normalsize Michigan Institute for Data Science \\
\normalsize Center for the Study of Complex Systems \\
\normalsize University of Michigan \\
\normalsize \texttt{encc@umich.edu}
\and
\textbf{Kevin Josey} \\
\normalsize Department of Biostatistics \\
\normalsize Harvard T.H. Chan School of Public Health \\
\normalsize \texttt{kjosey@hsph.harvard.edu} \\
\and
\textbf{Fan Yang} \\
\normalsize Department of Biostatistics and Informatics \\
\normalsize Colorado School of Public Health \\
\normalsize \texttt{fan.3.yang@cuanschutz.edu} \\
\and
\textbf{Debashis Ghosh}\\
\normalsize Department of Biostatistics and Informatics \\
\normalsize Colorado School of Public Health \\
\normalsize \texttt{debashis.ghosh@cuanschutz.edu}
}

\maketitle

\begin{abstract}
We propose a new perspective for the evaluation of matching procedures by considering the complexity of the function class they belong to. Under this perspective we provide theoretical guarantees on post-matching covariate balance through a finite sample concentration inequality. We apply this framework to coarsened exact matching as well as matching using the propensity score and suggest how to apply it to other algorithms. Simulation studies are used to evaluate the procedures.
\end{abstract}
\textbf{keywords:} {Causal effects, empirical distribution function, entropy metric, superpopulation, tail inequality, Vapnik-Chervonenkis dimension.}

\section{Introduction}
\label{sec:introduction}

\noindent Causal inference is a central goal for outcomes and policy research, particularly in the medical field. Among the many topics in this broad field of study are methods for evaluating treatment effects with non-randomized data. There is an abundance of observational data in nearly every discipline of science. However, bias induced by confounding is inherent in observational studies. In this context, the researcher must account for every potential confounder in some way before they can establish causality. While randomization remains the gold-standard for inference, as there is no confounding by definition, randomizing individuals into treatment groups is often cost prohibitive and sometimes unethical for certain study designs.

Under the potential outcomes framework \citep{neyman1923applications, rubin1974estimating}, \citet{rosenbaum1983central} were able to describe how the propensity score plays a key role in causal effect estimation and inference with observational data. The propensity score is defined as the probability of receiving a treatment given a set of measured covariates. Under strong ignorabiligy assumption, the propensity score removes bias attributable to confounding due to its property as a balancing score \citep{rosenbaum1983central}. 
With this result in mind, numerous methods for causal effect estimation were subsequently developed around the propensity score, with covariate balance serving as the primary objective (e.g., \cite{imai2014covariate,zubizarreta2015stable,chan2016globally}). However, the results presented by \citet{rosenbaum1983central} about the propensity score are derived in an asymptotic setting. This means that estimates of the propensity score may not adequately balance the covariate distribution in finite settings. {Therefore, many methods are resolved by iterating} between fitting a model for the propensity score and evaluating balance diagnostics on the propensity score adjusted covariates before estimating the treatment effect of interest. Some methods for evaluating  balance diagnostics have been proposed by \citet{ho2007matching} and \citet{sekhon2008multivariate}. The propensity score literature has mostly diverged into two overlapping yet distinct domains - one that uses the propensity score to derive balancing weights \citep{hainmueller2012entropy, imai2014covariate, chan2016globally}  and the other that uses a balancing score, such as the propensity score, to construct a matched cohort.


Recently, a multivariate matching approach using coarsened values of the observed covariates was developed by \citet{iacus2011multivariate}. They refer to their algorithm as coarsened exact matching. One of the primary aims of their method was to eliminate the iterative step of re-matching participants until an acceptable amount of balance is achieved. Coarsened exact matching is quite simple in nature and proceeds using the following high-level heuristic:
\begin{enumerate}
\item For each confounding variable, coarsen it into a certain number of categories;
\item Create strata based on the possible combinations of the coarsened values;
\item Compute a causal effect by comparing the outcomes of the treatment groups within the strata 
and adjusting for the stratum effect appropriately.
\end{enumerate}

The theoretical justification {provided by} \citet{iacus2011multivariate} for coarsened exact matching is a concept they term monotonic imbalance.  They show that bounding the distance between confounders to be small leads to matching procedures that are more flexible than procedures based on the equal percent bias reduction theory developed by Rubin and collaborators \citep{rubin1976multivariate, rubin1992affinely, rubin2006affinely}. One of the main advantages of coarsened exact matching is that it becomes amenable to large-scale database querying approaches to peforming causal inference: see \citet{salimi2016zaliql} as well as \citet{wang2017flame}.


However, fewer technical results exist for matching estimators than for other approaches, such as inverse probability weighting estimators. \citet{abadie2006large} have studied the large-sample asymptotics of matching estimators and found that in general, matching-based estimators of average causal effect did not have the usual $n^{1/2}$ convergence.  The intuition is that the matching algorithm introduces a bias into causal effect estimation that did not vanish asymptotically.  This bias term also increased with the number of confounders.  Bias-corrected estimators have been proposed by \citet{abadie2011bias}. \citet{abadie2016matching} performed a theoretical study of the asymptotic behavior of average causal effect estimators that match using the estimated propensity score.   

Conceptually, achieving covariate balance is a multivariate concept. If we let ${\cal L}(Z\mid T = 0)$ and ${\cal L}(Z\mid T = 1)$ denote the probability laws for the confounders conditional on treatment status then, ideally, as in the case of perfect randomization, these distributions are equal in some sense. We refer to this sense of equality as covariate balance.

Most covariate balance methods do not take the joint distribution of confounders into account but rather seek to match moments of the marginal distributions for the confounders. For example, \citet{imai2014covariate} proposed matching the first and second moments of covariates in their algorithm. Practically, one-dimensional diagnostics such as mean comparisons of confounders between treatment groups or Kolmogorov-Smirnov statistics are used to evaluate balance. \citet{wang2019minimal} have argued that due to the inherent complexity in attempting to achieve multivariate balance, one should instead strive to achieve approximate balance between confounders.

In this paper, we propose a new theoretical approach to evaluating and understanding covariate balance.
We introduce a distance metric to assess how close two multivariate distributions are from each other and define covariate balance as having zero distance.
This metric is defined in terms of the function family the matching procedure belongs to. Subsequent assessment of balance relies on understanding the behavior of the function classes in question. We demonstrate the following in the current paper:

\begin{enumerate}
\item The use of function classes fits naturally with the use of probability metrics \citep{zolotarev1984probability} for comparing probability laws and in this instance, multivariate distributions for confounders conditional on treatment.  
\item Results from empirical process theory \citep{van1996weak, kosorok2007introduction} can subsequently be used to study the behavior of function classes and to make probabilistic statements on the rates of convergence of matching procedures under ideal balance.
\item Ideal balance provides a new theoretical out-of-sample justification for the methodology of \citet{iacus2011multivariate} and can be used for the evaluation of other algorithmic strategies.  
\end{enumerate}
Based on the framework, one can view the techniques in this paper as being akin to developing a scalable strategy for achieving covariate balance that has relatively low complexity from the viewpoint described in Section \ref{sec:main}.

\section{Background and Preliminaries}
\label{sec:background}

\subsection{Data Structures and Causal Estimands}

Let the data be represented as $(Y_i,T_i,{Z}_i)$, $i=1,\ldots,n$, a
random sample from the triple $(Y,T,{Z})$, where $Y$
denotes the response of interest, $T$ denotes the treatment group, and $Z$ is a $p$-dimensional vector of covariates. We assume that
$T$ takes values {in} $\{0,1\}$.  

We now briefly review the potential 
outcomes framework \citep{rubin1974estimating, holland1986statistics}.  
Let $\{Y(0),Y(1)\}$ denote the potential outcomes for all $n$ subjects, and
the observed response {be} related to the potential outcomes by
\[Y = (1-T)Y(0) + TY(1).\]
In the potential outcomes framework, causal effects are defined as within-individual contrasts based on the potential outcomes.  
One popularly used estimand is the average causal effect, defined as 
\[
\text{ACE} = \frac{1}{n} \sum_{i=1}^n \left(Y_i(1) - Y_i(0)\right).
\]

Many assumptions are needed for performing valid causal inference.  These include the consistency assumption, the treatment positivity assumption, and the strongly ignorable treatment assumption \citep{rosenbaum1983central}, defined as 
\begin{equation}\label{unc}
T \perp \{Y(0),Y(1) \} \mid {Z}.
\end{equation}
Assumption (\ref{unc}) means that treatment assignment is conditionally independent of the set of potential outcomes given the covariates. Treatment positivity refers to $1 > P(T = 1\mid Z) > 0$ for all values of $Z$. Thus, the intuition is that any individual can potentially receive either treatment. Finally, the consistency assumption ensures that the observed outcome and the potential outcome under the observed treatment coincide.    

As described recently by \citet{imbens2015causal}, causal inference proceeds by modelling the assignment mechanism using observed covariates.   A quantity that naturally arises from this modelling is the propensity score \citep{rosenbaum1983central}, the probability of receiving treatment given confounders.  
The propensity score is defined as
\[
e({Z}) = P(T = 1\mid {Z}).
\]
Given the treatment ignorability
assumption in (\ref{unc}), it also follows by Theorem 3 of \citet{rosenbaum1983central} that treatment is
strongly ignorable given the propensity score, i.e.
\[ {T} \perp
\{Y(0),Y(1)\} \mid e(Z).\] 

Based on these assumptions and definitions, we can formulate causal inference using the following approach: (a) define an appropriate causal estimand; (b) formulate a propensity score model;
(c) check for covariate balance; (d) if (c) holds, estimate the causal estimand by conditioning on the propensity scores. We note that steps (b) and (c) tend to be iterative in practice. While the results in this paper pertain to propensity-matched analyses, they apply to more general matching strategies as well.
\subsection{Previous results on covariate balance}\label{sec:previous}

In terms of covariate balance, a major class of theoretical results come from work on {\it equal percent bias reduction} procedures (Rubin and Thomas, 1992, 1996).  Equal percent bias reduction means that a certain type of covariate matching will reduce bias in all dimensions of ${Z}$ by the same amount.

Define a matching method to be affinely invariant if the matching procedure is invariant to affine transformations of the covariates. If ${Z}$ given $T$ is assumed to have a so-called elliptically symmetric distribution, then Theorem 3.1. and Corollaries 3.1. and 3.2 of \citet{rubin1992affinely} apply so that any affinely invariant matching method will be equal percent bias reducing.   Examples of elliptically symmetric distributions include the multivariate normal and t distributions.   While elliptical symmetry of the confounders given treatment group is a restrictive assumption, this was relaxed in more recent work by \citet{rubin2006affinely}. There, they assumed that the conditional distribution of ${Z}$ given $T$ is a discriminant mixture of elliptically symmetric distributions. \citet{rubin2006affinely} prove that a generalization of equal percent bias reducing holds for this setup as well. 

Thus, for equal percent bias reducing methods, we have a guarantee that attempting to increase balance in one variable will not lead to distortions in balance for other variables.  However, the assumptions needed for equal percent bias reducing to hold seem restrictive in practice. \citet{iacus2011multivariate} took another approach by focusing on in-sample covariate discrepancies and requiring that the maximum discrepancy in sample means between treated and control subjects be bounded above by a constant.  They generalize this to arbitrary functions of the data, which they term imbalance bounding and define monotonic imbalance bounding matching methods to be those in which the discrepancies between a monotonic function applied to a variable is bounded above by a confounder-specific term.
Thus, one can be more stringent in the balance in variable without impacting the maximal imbalance across all confounders.  

There are many important implications of  requiring the monotonic imbalance bounding property.  First, many methods  of confounder adjustment, such as nearest-neighbor or caliper matching as defined in Cochran and Rubin (1973), are not monotonic imbalance bounding because they fix the number of treated and control observations within strata, while monotonic imbalance bounding methods imply variable numbers of observations.  By contrast, if the caliper matching procedure were to allow for different calipers for each confounder, then this would be monotonic imbalance bounding.  

\citet{iacus2011multivariate} also show that a key goal in causal effect estimation is to reduce model dependence \citep{ho2007matching}, meaning that there should not be  extrapolation of potential outcomes to regions in the covariate space where there are no observations. Under some assumptions on the model for potential outcomes, they show that for monotonic imbalance bounding methods, the model dependence is upper bounded by terms involving an imbalance parameter.
In addition, the estimation error for average causal effects using monotonic imbalance bounding matching methods can also be upper bounded by terms involving this parameter.

As a concrete example of a new monotonic imbalance bounding method, \citet{iacus2011multivariate} propose a coarsened exact matching algorithm for creating strata.  It proceeds as follows:
\begin{enumerate}
\item For each variable $Z_j$ $(j=1,\ldots,p)$, coarsen it into a function $C_j(Z_j)$ which takes on fewer values than the unique values of $Z_j$;
\item Perform exact matching between treated and control observations using the vector 
\[\left(C_1(Z_1),C_2(Z_2),\ldots,C_p(Z_p)\right).\]
This effectively creates strata ${\cal S}_1,\ldots,{\cal S}_J$ based on the unique combinations of 
\[\left(C_1(Z_1),C_2(Z_2),\ldots,C_p(Z_p)\right).\]
\item Discard strata in which there are only observations with $T = 0$.  For strata with only observations from the $T = 1$ population, extrapolate the potential outcome $Y(0)$ using the available controls or discard by restricting the causal effect of interest on the treated units for which causal effect can be identified without further modelling based assumptions. For strata with both treated and control observations, compare the outcome between the two populations.   
\end{enumerate}
\citet{iacus2011multivariate} have developed very easy-to-use software packages for implementing coarsened exact matching in R and Stata. They show that the coarsened exact matching approach satisfies the monotonic imbalance bounding property with respect to a variety of functionals of interest.  In addition, they provide a very intuitive explanation for what coarsened exact matching attempts to mimic. While classical propensity score approaches attempt to mimic a randomized study, analyses using coarsened exact matching will mimic randomized block designs, where the blocks are by definition predictive of the potential outcomes.  It is well-known that in this situation, randomized block designs will yield more efficient estimators (e.g., Box, Hunter and Hunter, 1978).

The other approach that has become of recent interest has been to incorporate covariate balance as part of the causal effect estimation process. For example, \citet{imai2014covariate} propose using  generalized methods of moments for causal effect estimation in which covariate balance is treated as a constraint in the procedure.  \citet{chan2016globally} propose the use of calibration estimators for causal effect estimation in which covariate balance constraints lead to a constrained Lagrangian dual optimization problem.  For these approaches, the authors are able to develop consistency and asymptotic normality results for the causal effect estimators. 

As described in more detail in  Section~\ref{sec:ideal_balance}, we will be using an integral probability metric to assess covariate balance among the two populations. In \citet{kallus2020generalized} a similar metric is used. They define such a metric as the target error to be minimized for obtaining optimal weighting coefficients when estimating the \textit{sample average treatment effect on the treated}. While our approaches are complementary, there are several notable differences. First, in \citet{kallus2020generalized}, they use their metric to find weights that correspond to known matching methods. The functions involved in their metric represent the expected relationship between potential outcomes and covariates. In our case, we take any matching procedure and given the measure of match, bound it by the probability metric involving functions representing the matching procedure itself, and provide probability bounds to how good the matching is. In addition, in  \citet{kallus2020generalized}, they assume a fixed population and therefore no randomness in covariate values, while our concern indeed focuses on the sample distribution of these covariates. The difference between these two approaches is further explained in Section~\ref{sec:modes_inference}.

\subsection{Modes of inference and covariate balance}
\label{sec:modes_inference}

In looking at the various proposals for accommodating covariate balance, it is useful to reconsider the ways in which one can perform causal inference.   \cite{imbens2015causal} have a nice overview on the distinction between finite-population and superpopulation modes for causal inference.   The finite-population mode of causal inference treats the sampled units as the population of interest.  The stochastic nature of the experiment is due solely to the treatment mechanism so that randomness occurs only with respect to the treatment assignments.   If one adopts the finite-sample point of view for causal inference, then one can use a randomization-based approach to performing inference for causal effects.  

By contrast, the superpopulation mode of inference considers two sources of variability.  The first is due to the randomness in the treatment assignments, and the second is due to the fact that the sampling units are a random sample from a superpopulation.  Thus, this approach posits a superpopulation from which the sampling units come from.  

Revisiting the previous work from \ref{sec:previous}, the equal percent bias reduction theory and the work of \cite{iacus2011multivariate} posit results about covariate balance assuming a finite-population mode for causal inference.  Thus, covariate balance results of these methods will involve subsampling and matching from the sampling units, and the balance occurs with respect to the matched sample.  The concept of balance we introduce in the next section can accommodate both modes of inference.


\section{Main Results}
\label{sec:main}

\subsection{Ideal Balance}
\label{sec:ideal_balance}
In this section, we wish to study covariate balance from the viewpoint of comparing the distributions ${\cal L}(Z\mid T = 0)$ and ${\cal L}(Z\mid T = 1)$. To do so, we must determine how this comparison is done. We do this by first defining probability pseudometrics.

\begin{definition}[Pseudometric]
Let $\mathcal{A}$ be the set of probability measures defined on a shared measurable space. A function $m:\mathcal{A}\times \mathcal{A} \rightarrow [0,\infty)$ is a {\bf pseudometric} on $\mathcal{A}$ if, for all $\mu$, $\nu$, $\lambda$ $\in$ $\mathcal{A}$, the following conditions are satisfied:
\begin{enumerate}
    \item $m(\mu,\mu) = 0$.
    \item $m(\mu,\nu) = m(\nu,\mu)$.
    \item $m(\mu,\nu)\leq m(\mu,\lambda) + m(\lambda,\nu)$.
\end{enumerate}
\end{definition}

Note these properties almost make $m$ a metric on $\mathcal{A}$, but notably we do not assume that if the distance between two elements is zero, then the two elements are the same. For the purpose of this paper, we will abuse terminology and refer to pseudometrics as metrics.

The class of metrics we will work with in this article is given by
\begin{equation}\label{ksfun}
\gamma_{\mathcal{F}}(\mu,\nu) = \sup_{f \in {\cal F}} \left| \int f d\mu - \int f d\nu \right|,
\end{equation}
where ${\cal F}$ is a class of functions. In (\ref{ksfun}), $\gamma_{\mathcal{F}}(\mu,\nu)$ is referred to by \citet{zolotarev1984probability} as an example of a probability metric. In our notation, we drop the dependency of $\gamma_{\mathcal{F}}$ on $\mathcal{F}$ and write it as $\gamma$. We now define ideal balance as being based on (\ref{ksfun}).

\begin{definition}[Ideal Balance]\label{def:idealbalance}
Let $\mu$ and $\nu$ be distributions on the same probability space and $m$ a pseudometric, then we say $\mu$ and $\nu$ satisfy {\bf Ideal Balance} with respect to $m$ if $m(\mu,\nu) = 0$.
\end{definition}

{When $\mu$ and $\nu$ are the conditional distributions of the covariates given the treatment group, as in Section~\ref{sec:background}, ideal balance is a restriction on the population. If these are instead the empirical distributions of the data, ideal balance is a sample restriction. Matching methods, in a sense, intend to achieve ideal balance on the matched data for some $m$.}

Note that at this stage, we have only dealt with population distributional laws and have not described how to estimate or compute these quantities with real data. In practice, we would not expect ideal balance to hold in observational studies. However, it does serve as a useful benchmark through which we can study the behavior of various functional constraints. Here, the function spaces ${\cal F}$ in (\ref{ksfun}) play the role of the constraints; {more complex function spaces correspond to more constraints on the joint distributions} of $Z|T = 1$ and $Z|T = 0$.

\subsection{A Concentration Inequality Result}

Let $\mathcal{F}$ be a function space and $\| \cdot \|$ a norm. The covering number $N(\epsilon,{\cal F},\| \cdot \|)$ is the minimum number of $\|\cdot\|$-balls of radius $\epsilon$ needed to cover ${\cal F}$, where a ball centered around $f\in \mathcal{F}$ is the set $\{g \mid \|f - g\| \leq \epsilon\}$. Intuitively, one can think of the covering number as a measure of the complexity of the function class ${\cal F}$. For a measure $\mu$ the norm $L_r(\mu)$-norm, for $r\geq 1$, is defined as $\|f \|_{L_r(\mu)}^r = \int |f|^r d\mu$. Throughout the paper, we will assume $\mathcal{F}$ is uniformly bounded. Note that if $\mu$ is any probability measure, and under uniform boundedness, we can endow $\mathcal{F}$ with the norm $L_r(\mu)$ without dropping any of its elements. Unless otherwise specified, we assume the range of the functions in ${\cal F}$ is $[0,1]$. Finally, for a function class ${\cal F}$, an envelope function of $\mathcal{F}$ is defined as any function $h$ such that for all $f$ in $\mathcal{F}$, the inequality
\[
|f(x)| \leq |h(x)|
\]
is satisfied for any $x$.

Let $\{Z_i\}_{i=1}^n$ be a sample where each $Z_i$ has distribution $Q$. We denote the empirical distribution by $\mathbb{Q}_n$. The $\mathcal{F}$-indexed empirical process $\mathbb{G}_n^Q$ is defined as the map taking any $f\in\mathcal{F}$ to
$$
\mathbb{G}_n^Q(f) = \sqrt{n}\left(\int fd\mathbb{Q}_n - \int fdQ\right) =\frac{1}{\sqrt{n}} \sum_{i=1}^n \left(f(Z_i) - \int f dQ\right).
$$

\begin{theorem}
\label{thm:concentration_theorem}
Let $\mathbb{Q}^0_{n_0}$ and $\mathbb{Q}^1_{n_1}$ be two empirical distributions of observations sampled from $Q^0$ and $Q^1$, respectively, and assume ideal balance holds for $Q^0$ and $Q^1$ with respect to $\gamma$. Let $M$ be the collection of probability measures.  If there exists constants $C$ and $K$ such that  ${\cal F}$ satisfies
\[ \sup_{{\mu\in M}}  N(\epsilon,{\cal F},\|\cdot\|_{L_r({\mu})} ) \leq \left ( \frac{K}{\epsilon} \right)^C,\]
for every $0 < \epsilon < C$, then
\begin{equation}\label{concineq}
Pr \{ \gamma(\mathbb{Q}^0_{n_0},\mathbb{Q}^1_{n_1}) > \delta\} \leq \left ( \frac{D\delta}{2\sqrt{C}} \right)^{C} \left( n_0^{C/2}\exp(-n_0\delta^2/2) + n_1^{C/2}\exp(-n_1\delta^2/2) \right),
\end{equation}
where $D$ is a constant depending on $K$ only.
\end{theorem}
The proofs of Theorem~\ref{thm:concentration_theorem} and subsequent results are found in the supplementary material. Throughout the paper, we will use $B_n(\delta,D,C)$ for the bound in Theorem~\ref{thm:concentration_theorem}, where the subscript $n$ reminds us of the dependence on the sample size.

\begin{remark}
We note that the bound in (\ref{concineq}) is nonasymptotic and will hold for any sample size.
\end{remark}

\begin{remark}
In this framework, the function classes play an important role. Theorem~\ref{thm:concentration_theorem} gives a bound in terms of the entropy number of the function class in question.
In particular, low-complexity functions are favored using this approach.  A key technical point is ensuring that the covering number condition in the theorem is satisfied.  To do so, we will primarily use results from Vapnik-Chervonenkis theory \citep{chervonenkis1971uniform} to determine appropriate covering numbers.
\end{remark}

In most cases the function classes of interest are not real-valued but vector-valued. The following straightforward results can be used to deal with these cases.

\begin{lemma}\label{lemma:gammasum}
Let $\{\mathcal{F}_i \}_{i=1}^d$ be a collection of real-valued function spaces and $({P}^{i},{Q}^{i})$ satisfy ideal balance under $\gamma_{\mathcal{F}_i}$ for each $1\leq i \leq d$. Let $(\mathbb{P}^i,\mathbb{Q}^i)$ denote their respective empirical distributions with implicit sample size dependence. Then
\[
    Pr\left( \sum_{i=1}^d \gamma_{\mathcal{F}_i}(\mathbb{P}^i,\mathbb{Q}^i) > \delta\right) \leq \sum_{i=1}^d B(\delta/d,D_i,C_i).
\]
\end{lemma}

Now, consider the collection $\{\mathcal{F}_i \}_{i=1}^d$, where each $\mathcal{F}_i$ is a real-valued function space. Define $\mathcal{F}=\{f = (f_1, \dots,f_d)^T \mid f_i\in\mathcal{F}_i \,\,\, \text{for all} \,\,\, i \}$. Let $\pi_\ell$ be the $\ell^{th}$ coordinate projection, that is, for a finite dimensional vector $x = (x_1, \dots, x_d)$, $\pi_\ell(x) = x_\ell$. Finally, define $\mathcal{F_\pi}=\{\pi_\ell\circ f \mid f\in\mathcal{F}, 1\leq\ell\leq d \}$. Note the elements of $\mathcal{F}_\pi$ are real-valued. The following lemma tells us we can either assume $\mu$ and $\nu$ satisfy ideal balance with respect to each of $\gamma_{\mathcal{F}_i}$, or that they satisfy ideal balance with respect to $\gamma_{\mathcal{F}_\pi}$.

\begin{lemma}\label{lemma:fequiv}
Let $\mathcal{F}$, $\{\mathcal{F}_i\}_{i=1}^d$, and $\mathcal{F}_\pi$ be as above, and let $\mu$ and $\nu$ denote two probability measures.  Then the following are equivalent:
\begin{enumerate}
    \item $\mu$ and $\nu$ satisfy ideal balance with respect to $\gamma_{\mathcal{F}_\pi}$;
    \item $\mu$ and $\nu$ satisfy ideal balance with respect to each $\gamma_{\mathcal{F}_i}$, $1\leq i\leq d$.
    \item $\max_{i}\gamma_{\mathcal{F}_i}(\nu,\mu) = 0$.
\end{enumerate}
\end{lemma}

The following corollary will be very useful:
\begin{corollary}\label{corollary:lpnorm}
Let $\mathcal{F}$ and $\mathcal{F}_\pi$ be as above, and $\mathcal{F}_i= \mathcal{F}^*$ for all $i$. Assume $\mathcal{F}^*$ has polynomial covering number. Let $\{X^0_j\}_{j=1}^{n_0}\sim Q^0$ and $\{X^1_j\}_{j=1}^{n_1} \sim Q^1$, where $Q^0$ and $Q^1$ satisfy ideal balance with respect to $\gamma_{\mathcal{F}_\pi}$. Fix $f^*\in\mathcal{F}$, then
\[
    Pr\left( \left\Vert \frac{1}{n_0}\sum_{j=1}^{n_0} f^*(X^0_j) - \frac{1}{n_1}\sum_{j=1}^{n_1} f^*(X^1_j) \right\Vert_{\ell_p}>\delta\right) \leq dB(\delta/{d^{1/p}},D^*,C^*),
\]
for finite $p\geq 1$, and
\[
    Pr\left( \left\Vert \frac{1}{n_0}\sum_{j=1}^{n_0} f^*(X^0_j) - \frac{1}{n_1}\sum_{j=1}^{n_1} f^*(X^1_j) \right\Vert_{\ell_\infty}>\delta\right) \leq dB(\delta,D^*,C^*),
\]
where $D^*,C^*$ depend only on $\mathcal{F}^*$.
\end{corollary}

\begin{definition}[Vapnik-Chervonenkis Dimension]
The {\bf Vapnik-Chervonenkis} dimension of a function class ${\cal F}$ on an ambient set ${\cal X}$  is the cardinality of the largest subset shattered by ${\cal F}$. A function class $\mathcal{F}$ shatters a set $S\in\mathcal{X}$ if for each possible $0-1$ labeling of the elements of $S$ there is at least one function $f\in\mathcal{F}$ that realizes such labeling.
\end{definition}

A key result we will use is an application of Theorem 2.6.7 of \citet{van1996weak}, which implies that if a function class ${\cal G}$ has finite Vapnik-Chervonenkis dimension $v$, then 
\[ \sup_{\mu}  N(\epsilon,{\cal G},L_2(\mu)) \leq \left ( \frac{K}{\epsilon} \right)^{C^*},\]
where $C^* = 2v - 2$.

\section{Examples}
\label{sec:examples}

\subsection{Balance on coarsened function classes}
Consider coarsened exact matching as described in \cite{iacus2011multivariate}. Let $\mathcal{Z}_0=\{Z^0_i\}_{i=1}^{n_0}$ and $\mathcal{Z}_1=\{Z^1_j\}_{j=1}^{n_1}$ be the control and treatment samples, respectively. In coarsened exact matching we create a partition of the sample space and match samples which are found in the same element of the partition, and discard samples in subsets without samples from the opposite group. We are interested in the quantity
$$
\Delta = \frac{1}{m_0}\sum_{i\in M_0} w_i^0 Z_i^0 - \frac{1}{m_1}\sum_{j\in M_1} w_j^1 Z_j^1,
$$
where $m_\ell$ is the number of matched samples for the $\ell^{th}$ group, $M_\ell$ is its index set, and $\{w_i^0,w_j^1\}_{i\in M_0, j\in M_1}$ are weights.

In the supplementary material we describe how to express this matching procedure as a function $f$ on the variables ${Z_i^0}$ and $Z_j^1$. This allows us to express $\Delta$ in terms of $f$. We further specify the function space $\mathcal{F}$ for which
$$
\|\Delta\|\leq\gamma_{\mathcal{F}}(\mathbb{Q}^0_{n_0},\mathbb{Q}^1_{n_1})
$$
holds for an appropriate norm.
Using the properties of $\mathcal{F}$ and provided the bound above, we can derive our results of interest:
$$
Pr(\left|\Delta_k\right|\geq\delta)\leq B(\delta,D,C^*),
$$
for a constant $C^*$ and where $\Delta_k$ is the $k^{th}$ component of $\Delta$. Similarly,
$$
Pr(\|\Delta\|_{\ell_p}\geq\delta)\leq dB(\delta/d^{1/p},D,C^*)
$$
and 
$$
Pr(\|\Delta\|_{\ell_\infty}\geq\delta)\leq dB(\delta,D,C^*).
$$

\subsection{Covariate balance on the linear propensity score}
\label{sec:linear_propensity}

As discussed in Section \ref{sec:main}, there has been a lot of work on developing matching results based on linear discriminant analysis. That is, we assume that $P(Z \mid T = \ell)$ follows $N(\mu_\ell, \Sigma)$. Under this model, the metric for consideration is the $logit$ of the propensity score (see \cite{stuart2010matching}). In the supplementary material we show the distance $\left| logit(e(Z)) - logit(e(Z')\right|$ can be expressed in terms of the linear discriminant analysis hyperplance vector. Indeed, if $p$ is the dimension of the covariates, we can create a function space $\mathcal{F}$ derived from hyperplanes and with Vapnik-Chervonenkis dimension $p+1$ such that
\begin{align*}
    \Delta &=
    \left| \frac{1}{m_0}\sum_{i\in M_0} logit(e(Z_i)) - \frac{1}{m_1}\sum_{j\in M_1} logit(e(Z_j)) \right| \\
    & \leq
    \gamma_{\mathcal{F}}(\mathbb{Q}^0_{n_0},\mathbb{Q}^1_{n_1}),
\end{align*}
allowing us, using Theorem~\ref{thm:concentration_theorem}, to determine the bound of interest:
\begin{equation*}
    Pr \{ \Delta > \delta\} \leq B(\delta,D,2p).
\end{equation*}


\subsection{Covariate balance using kernels}
Many authors \citep{hazlett2016kernel, wong2018kernel, zhu2018kernel} have advocated for the use of kernel methods for matching and evaluating covariate balance. This corresponds to assuming that ${\cal F}$ in (\ref{ksfun}) represents a Reproducing Kernel Hilbert space. Further details about these function spaces can be found in the supplementary material.

To apply Theorem~\ref{thm:concentration_theorem} to the kernel setting, we will note there exists a version of linear discriminant analysis from section~\ref{sec:linear_propensity} that can be extended to the reproducing Kernel Hilbert Space setting \citep{baudat2000generalized}. Let $\mathcal{H}$ be a reproductive kernel Hilbert space and $\|\cdot\|_{\mathcal{H}}$ the norm associated to it, then a natural metric to consider for a kernelized matching procedure would be 
\[
\Delta_{\mathcal{H}} = \left\Vert \frac{1}{m_0} \sum_{i \in M_0} f(Z_i) - \frac{1}{m_1} \sum_{j \in M_1} f(Z_j) \right\Vert_{\mathcal{H}},
\]
which represents a functional generalization of $\Delta$ from Section~\ref{sec:linear_propensity}, and where $f\in\mathcal{H}$ is an appropriate function chosen by the user. Then $\Delta_{\mathcal{H}} \leq \gamma_{\cal F}(\mathbb{Q}^0_{n_0},\mathbb{Q}^1_{n_1})$, and we can use the previous results with a few adjustments. We show in the supplementary material that
\[
P(\Delta_{\mathcal{H}} > \delta) \leq B(\delta,D,C^*),
\]
where $C^*$ depends on the smoothness properties of $\mathcal{H}$.

\section{Practical implementation}
\label{sec:practical}

So far, we have given theoretical results that describe how algorithms under various function classes behave under the ideal balance assumption. As noted earlier, the ideal balance definition is strict but permits theoretical characterization of various algorithms. The question then naturally arises as to how to use the theoretical results from the previous sections in practice.

Note one can view the metric in equation~(\ref{ksfun}) as a multivariate balance metric, which differentiates it from many other balance metrics in the literature. \citet{zhu2018kernel} used (\ref{ksfun}), where ${\cal F}$ is a reproducing kernel Hilbert space, as a covariate balance diagnostic. There, they found that in certain situations, the diagnostic was more sensitive in finding covariate imbalances relative to univariate diagnostics as well as those based on the prognostic score \citep{hansen2008prognostic}. 

Consider the problem of estimating the average causal effect among the treated. In practice, it is unlikely that ideal balance will hold for the treatment and control populations. That is to say, $\gamma_{\mathcal{F}}\left(Q^{0}, Q^{1}\right) \ne 0$, unless treatment is randomized. Therefore, we wouldn't be able to use Theorem~\ref{thm:concentration_theorem} in an observational study. However, a slight modification can be done for which the analysis remains largely the same.

Let $w \in \mathcal{W} \subset \mathbb{R}^{n_0}$ be a weight vector and define
\[
\mathbb{Q}^0_{w} = \frac{1}{\sum_{i:T_i = 0}w_i}\sum_{i:T_i = 0} w_i\delta_{X_i}.
\]
The majority of methods in causal inference have as a goal to find appropriate weights $w$ for which $\mathbb{Q}^0_{w}$ converges to $Q^*$ for some distribution $Q^*$ that indeed satisfies ideal balance with $Q^1$. That is, for which $\gamma_{\mathcal{F}}\left(Q^*, Q^1\right) = 0$. In order for this modification to be feasible, we just need to modify our proof of Theorem~\ref{thm:concentration_theorem} and include the convergence rates of $\mathbb{Q}^0_w$ to $Q^*$, which may change depending on the problem. Having done so, we continue in a parallel manner.

Let $f^*\in\mathcal{F}$ represent a matching procedure with balance diagnostic
\[
\Delta=\left| \int f d\mathbb{Q}_{w}^0 - \int f d\mathbb{Q}^1_{n_1} \right|,
\]
then, by the definition of $\gamma_\mathcal{F}$,
\[
\Delta \leq \gamma_{\mathcal{F}}\left(\mathbb{Q}^0_{w},\mathbb{Q}^1_{n_1}\right).
\]
Therefore, if we can find weights for which $\mathbb{Q}^0_{w}$ converges to $Q^*$ and $\gamma_{\mathcal{F}}(Q^*, Q^1) = 0$, then we can bound the probability that $\Delta$ exceeds some threshold $\delta$.

There are many methods for finding $w \in \mathcal{W}$, the most straightforward being the inverse probability of treatment weights,
\[ w_i = T_i + \frac{e(Z_i)(1 - T_i)}{1 - e(Z_i)}.\]
Even heavily prescribed matching algorithms that are found throughout the causal inference literature find some weights $w \in \mathcal{W}$ as described by \cite{abadie2006large}. In one-to-one matching with replacement, let $\mathcal{J}(i) = \{j_1(i), j_2(i), \ldots\}$ be the set of indices of units that are matched with the unit $i = 1,2,\ldots,n$. If there are no ties, then $\mathcal{J}(i) = j(i)$. With ties present, which occur frequently especially with exact matching (see coarsened exact matching), $\mathcal{J}(i)$ might contain multiple matched indices. The matching process will allow us to produce weights for every unit by solving 
\[
w_i = \sum_{\{l:T_l = 1\}} \frac{I[i \in \mathcal{J}(l)]}{\#\mathcal{J}(l)} \ \text{for all} \ i \in \{i:T_i = 0\}
\]
where $\#\mathcal{J}(i)$ denotes the cardinality of $\mathcal{J}(i)$.

\section{Simulation Studies}
\label{sec:simulations}

\begin{figure}[!ht]
	\centering
	\includegraphics[scale = 0.55]{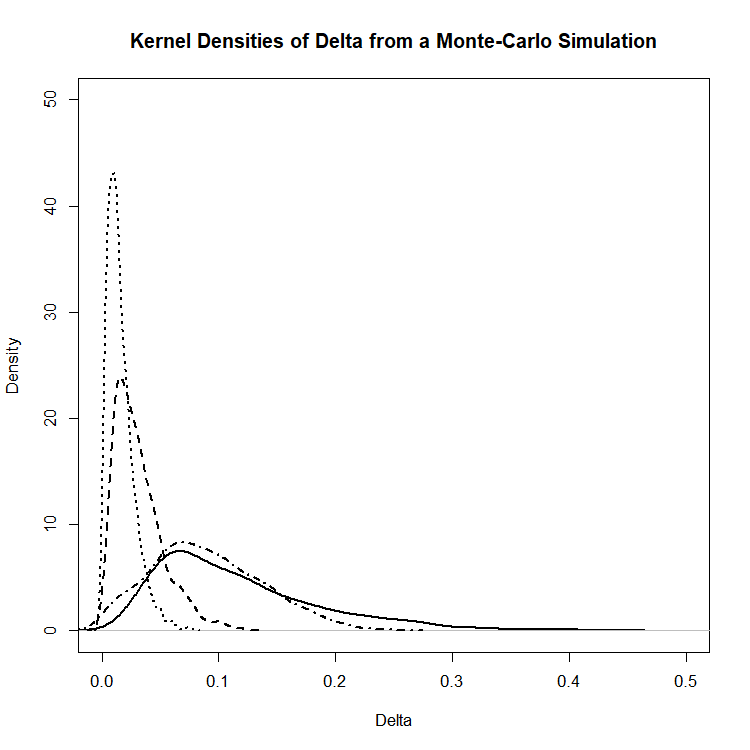}
	\caption{Kernel Densities of the $\Delta$ balancing statistics for the baseline scenario with $n = 1000$ and $\sigma^2 = 10$.  The solid line is the distribution from the unweighted estimates, the dashed line is the distribution for coarsened exact matching, the dotted line is the distribution for the linear propensity score, and the dotted-dashed line for the support vector machine examples.}\label{figure:kernel}
\end{figure}

We perform a simulation study to evaluate the distribution of the distances reported in Section \ref{sec:examples}. We also examine their downstream consequences for estimating average treatment effects on the treated. There are two data generating mechanisms that we consider. In addition, we vary the sample size and the variance of the responses for a total of eight scenarios. We replicate each of these scenarios, described below, over $1000$ iterations. We report the mean and Monte Carlo standard errors of the three distances ($\Delta$) examined in Section \ref{sec:examples} (Table \ref{table:delta}) along with the kernel density estimates for one representative scenario (Figure \ref{figure:kernel}). We also evaluate the downstream effects of these $\Delta$ statistics on the average treatment effect using one-to-one matching methods described by \citet{abadie2006large} implemented in the \texttt{Matching} package \citep{sekhon2008multivariate} (Tables \ref{table:estimate} and \ref{table:coverage}).

For $i = 1,2,\ldots,n$, let $Z_{i1} \sim \mathcal{N}(1, 4)$,  $Z_{i2} \sim \text{Bin}(1, 0.3)$, $Z_{i3} \sim \mathcal{N}(0, 1)$, and $Z_{i4} \sim \text{Bin}(1, 0.5)$ where $T_i$ denotes the binary treatment assignment. The conditional means of the outcomes for the treated, $\mu_1(Z_i)$, and the controls, $\mu_0(Z_i)$, are constructed as
\begin{equation}\label{mu}
\begin{split}
\mu_0(Z_i) &= 10 - 3Z_{i1} - Z_{i2} + Z_{i3} + 3Z_{i4} \enskip \text{and} \\
\mu_1(Z_i) &= \mu_0(Z_i) + 5 + 3Z_{i1} - Z_{i2} + Z_{i3} - 3Z_{i4}.
\end{split}
\end{equation}
We sample $T_i \sim \text{Bin}(1, 0.5)$ distribution. For $i = 1,2,\ldots,n$, we sample the counterfactual responses $Y_i(1) \sim \mathcal{N}[\mu_1(Z_i), \sigma^2]$ and $Y_i(0) \sim \mathcal{N}[\mu_0(Z_i), \sigma^2]$. The observed outcome is $Y_i = T_iY_i(1) + (1 - T_i)Y_i(0)$. We will refer to these conditions with the label ``baseline". For the error variance, we set $\sigma^2 \in \{5,10\}$. 

For the scenario labeled ``sparse", we include an additional set of covariates that ultimately do not affect the outcome. The outcomes are determined by the potential outcome models in (\ref{mu}), yet the methods we consider also account for the noise covariates $Z_{i5} \sim \mathcal{N}(-1, 4)$,  $Z_{i6} \sim \text{Bin}(1, 0.7)$, $Z_{i7} \sim \mathcal{N}(0, 1)$, and $Z_{i8} \sim \text{Bin}(1, 0.5)$.

As mentioned before, we test the three examples described in Section \ref{sec:examples} in their ability to produce efficient, unbiased estimates of the average treatment effect of the treated. Linear discriminant analysis sets $f$ to be the logit transformation of the fitted posterior probability that each unit receives treatment. The support vector machine examples use the distance that each point is from the resulting separating hyperplane assuming a linear kernel. Coarsened exact matching is performed {similar to what} is described in \cite{iacus2011multivariate} and is implemented with the \texttt{cem} R package. Table 1 shows the results of our simulation experiment. Since balance is already achieved through randomization in this simulation, we also report the unmatched, crude estimate of the average causal effect for references. Here the value $\Delta$ is the maximum absolute sample mean difference for the unweighted covariates.

\begin{table}[htbp!]
\centering
\begin{tabular}{cccccccc}
\hline
$n$ & $\sigma^2$ & Scenario & $\theta$ & A & B & C & D \\ \hline
1000 & 5 & baseline & 6.2 & 0.11 (0.07) & 0.03 (0.02) & 0.02 (0.01) & 0.09 (0.04) \\
1000 & 5 & sparse & 6.2 & 0.15 (0.07) & 0.01 (0.01) & 0.03 (0.02) & 0.13 (0.05) \\
1000 & 10 & baseline & 6.2 & 0.12 (0.07) & 0.03 (0.02) & 0.02 (0.01) & 0.09 (0.05) \\
1000 & 10 & sparse & 6.2 & 0.15 (0.07) & 0.01 (0.01) & 0.03 (0.02) & 0.13 (0.05) \\
2000 & 5 & baseline & 6.2 & 0.08 (0.05) & 0.02 (0.01) & 0.01 (0.01) & 0.06 (0.03) \\
2000 & 5 & sparse & 6.2 & 0.11 (0.05) & 0.01 (0.01) & 0.02 (0.01) & 0.09 (0.04) \\
2000 & 10 & baseline & 6.2 & 0.08 (0.05) & 0.02 (0.01) & 0.01 (0.01) & 0.06 (0.03) \\
2000 & 10 & sparse & 6.2 & 0.11 (0.05) & 0.01 (0.01) & 0.02 (0.01) & 0.09 (0.04) \\ \hline
\end{tabular}
\caption{Average and Monte Carlo standard error of $\Delta$ found in the experiment. In this table, Method A is the unweighted estimate, Method B refers to coarsened exact matching, Method C to linear discriminant analysis, and Method D to support vector machines. Since both A and B create a vector valued $\Delta$ we report the maximum.}\label{table:delta}
\end{table}

The values $\Delta$ are not necessarily directly comparable in this example. They do represent the distributions whose tail probabilities we are bounding in theorem. The simulation serves to characterize some of the densities of these statistics so that we might better understand which values of $\delta$ are acceptable for the different balance methods in Section \ref{sec:examples}. We see that the values for $\Delta$ after coarsened exact matching were the most heavily concentrated, followed closely by the values generated by linear discriminant analysis. The balance diagnostics from a support vector machine and from an unweighted comparison yielded considerably more dispersed values.

One point of direct comparison that we may take between the different $\Delta$ estimates is the downstream effects of the various balancing methods with estimating the average treatment effect. The purpose of this portion of the simulation study shows how the concentration of the distribution for $\Delta$ may have little to do with the actual quality of the average treatment effect estimates - the ultimate result for causal inference. Although the concentration of the distribution for $\Delta$ under coarsened exact matching was the most narrow among the other densities found for $\Delta$ under linear discriminant analysis and support vector machines, the estimated average treatment effect is also the most biased. The Monte Carlo standard errors also seem to be greater than the other two balance methods. Linear discriminant analysis also conferred a narrow concentration of $\Delta$ statistics yet produced the most efficient estimates of the average treatment effect, other than from the unweighted estimate which had the smallest Monte Carlo standard errors. This result is interesting because the unweighted diagnostics had the most dispersed values for $\Delta$. This leads us to believe that the scale of the $\Delta$ statistics must be carefully considered while evaluating balance to make some determination on which method is most suitable for evaluating treatment effects.

\begin{table}[htbp!]
\centering
\begin{tabular}{cccccccc}
\hline
$n$ & $\sigma^2$ & Scenario & $\theta$ & A & B & C & D \\ \hline
1000 & 5 & baseline & 6.2 & 6.20 (0.33) & 6.24 (0.33) & 6.20 (0.42) & 6.20 (0.36) \\
1000 & 5 & sparse & 6.2 & 6.20 (0.34) & 6.29 (1.24) & 6.21 (0.45) & 6.20 (0.39) \\
1000 & 10 & baseline & 6.2 & 6.20 (0.37) & 6.22 (0.40) & 6.20 (0.47) & 6.20 (0.42) \\
1000 & 10 & sparse & 6.2 & 6.19 (0.35) & 6.31 (1.46) & 6.20 (0.46) & 6.22 (0.42) \\
2000 & 5 & baseline & 6.2 & 6.19 (0.24) & 6.21 (0.24) & 6.20 (0.29) & 6.20 (0.25) \\
2000 & 5 & sparse & 6.2 & 6.20 (0.23) & 6.34 (0.71) & 6.21 (0.29) & 6.21 (0.26) \\
2000 & 10 & baseline & 6.2 & 6.21 (0.25) & 6.21 (0.26) & 6.19 (0.32) & 6.21 (0.28) \\
2000 & 10 & sparse & 6.2 & 6.21 (0.25) & 6.38 (0.79) & 6.21 (0.31) & 6.21 (0.27) \\ \hline
\end{tabular}
\caption{Summary of simulation estimates and Monte Carlo standard errors. The simulation scenarios corresponding to "baseline" and "sparse" are described in further detail in Section \ref{sec:simulations}. Here, $\theta$ refers to the population average treatment effect among the treated.  In this table, Method A is the unweighted estimate, Method B refers to coarsened exact matching, Method C is linear discriminant analysis, and Method D is support vector machines.}\label{table:estimate}
\end{table}

\begin{table}[htbp!]
\centering
\begin{tabular}{cccccccc}
\hline
$n$ & $\sigma^2$ & Scenario & $\theta$ & A & B & C & D \\ \hline
1000 & 5 & baseline & 6.2 & 0.952 & 0.937 & 0.941 & 0.929 \\
1000 & 5 & sparse & 6.2 & 0.944 & 0.955 & 0.934 & 0.917 \\
1000 & 10 & baseline & 6.2 & 0.941 & 0.918 & 0.935 & 0.912 \\
1000 & 10 & sparse & 6.2 & 0.955 & 0.950 & 0.951 & 0.931 \\
2000 & 5 & baseline & 6.2 & 0.931 & 0.945 & 0.937 & 0.923 \\
2000 & 5 & sparse & 6.2 & 0.956 & 0.945 & 0.939 & 0.918 \\
2000 & 10 & baseline & 6.2 & 0.959 & 0.936 & 0.926 & 0.928 \\
2000 & 10 & sparse & 6.2 & 0.953 & 0.946 & 0.948 & 0.935 \\ \hline
\end{tabular}
\caption{Summary of coverage probabilities from the simulation experiment. The simulation scenarios corresponding to "baseline", "interaction", "positivity", and "sparse" are described in further detail in Section \ref{sec:simulations}.   Here, $\theta$ refers to the population average treatment effect among the treated.  In this table, Method A is the unweighted estimate, Method B refers to coarsened exact matching, Method C to linear discriminant analysis, and Method D to support vector machines. }
\end{table}\label{table:coverage}

\section*{Acknowledgments}
The authors would like to acknowledge funding support from the following sources: the National Institutes of Health, the National Science Foundation, the Veterans Administration and the Grohne-Stepp Endowment from the University of Colorado Cancer Center.  

\section*{Appendix}

\subsection*{Proof of theorem \ref{thm:concentration_theorem}}
We will use $P$ and $Q$ instead of $Q^0$ and $Q^1$ to ease symbolic burden on the reader.
\begin{proof}
By definition of $\gamma$:
\begin{eqnarray*}
\gamma(\mathbb{P}_{n_0},\mathbb{Q}_{n_1}) &=& \sup_{f \in {\cal F}} \left| \int f d\mathbb{P}_{n_0} - \int f d\mathbb{Q}_{n_1} \right| \\
&=& \sup_{f \in {\cal F}} \left| \int f d \mathbb{P}_{n_0} \pm \int f dP \pm \int{fdQ} - \int f d\mathbb{Q}_{n_1} \right| \\
&\leq& \sup_{f \in {\cal F}}
\left|
\int f d \mathbb{P}_{n_0} - \int f dP - \int f d\mathbb{Q}_{n_1} + \int{fdQ}
\right|
+ \sup_{f \in {\cal F}}
\left|
\int f dP - \int{fdQ}
\right| \\
&=&  \sup_{f \in {\cal F}}
\left|
\int f d \mathbb{P}_{n_0} - \int f dP - \int f d\mathbb{Q}_{n_1} + \int{fdQ}
\right|,
\end{eqnarray*}
since $\gamma(P,Q) = 0$. Using elementary probability arguments, we have
\begin{eqnarray}\label{sumprobs}
Pr \{ \gamma(\mathbb{P}_{n_0},\mathbb{Q}_{n_1}) > \delta\} &=& Pr \left( \sup_{f \in {\cal F}}
\left|
\int f d \mathbb{P}_{n_0} - \int f dP - \int f d\mathbb{Q}_{n_1} + \int{fdQ}
\right| > \delta \right) \nonumber \\
&=& Pr \left( \sup_{f \in {\cal F}} \left| \frac{1}{\sqrt{n_0}}\mathbb{G}_{n_0}^P(f) - \frac{1}{\sqrt{n_1}}\mathbb{G}_{n_1}^Q(f) \right| > \delta\right) \nonumber \\
&\leq & Pr \left(\sup_{f \in {\cal F}} | \mathbb{G}_{n_0}^P(f) | > \sqrt{n_0}\delta/2\right) + Pr \left(\sup_{f \in {\cal F}} | \mathbb{G}_{n_1}^Q(f) | > \sqrt{n_1}\delta/2\right), \nonumber
\end{eqnarray}
where $ \mathbb{G}_{n_0}^P(f)$ and $ \mathbb{G}_{n_1}^Q(f) $ represent the $\mathcal{F}$-indexed empirical processes of $P$ and $Q$, respectively.

Applying Theorem 2.14.9 in \citet{van1996weak}, we can bound each of the terms as follows:
\[ Pr \left(\sup_{f \in {\cal F}} | \mathbb{G}_{n_0}^P(f) | >\sqrt{n_0}\delta/2 \right) <  \left ( \frac{D\sqrt{n_0}\delta}{2\sqrt{C}} \right)^C \exp(-n_0\delta^2/2)\]
\[ Pr \left(\sup_{f \in {\cal F}} \left| \mathbb{G}_{n_1}^Q(f) \right| > \sqrt{n_1}\delta/2\right) < \left( \frac{D\sqrt{n_1}\delta}{2\sqrt{C}}\right)^C \exp(-n_1\delta^2/2),\]
where $D$ is a constant depending only on $K$.  Plugging these two bounds into (\ref{sumprobs}) concludes the proof.
\end{proof}

\subsection*{Proof of Lemma \ref{lemma:gammasum}}

\begin{proof}
Define $\gamma_i=\gamma_{\mathcal{F}_i}(\mathbb{P}^i,\mathbb{Q}^i)$. Then:
\begin{align*}
    Pr\left(\sum_i\gamma_i>\delta\right) &= 1 - Pr\left(\sum_i\gamma_i<\delta\right) \\
    & \leq 1 - Pr(\gamma_i<\delta/d \,\,\, \forall i) \\
    & = Pr(\exists \,\,\,i \ni \gamma_i > \delta/d) \\
    & \leq \sum_i Pr(\gamma_i > \delta/d) \\
    & \leq \sum_i B(\delta/d,D_i,C_i),
\end{align*}
where we have used the union bound in the second inequality.
\end{proof}

\subsection*{Proof of Lemma \ref{lemma:fequiv}}
\begin{proof}
Assume $\gamma_{\mathcal{F}_i}(\mu,\nu) = 0$ for all $i$. Then
\begin{align*}
    \gamma_{\mathcal{F}_\pi}(\mu,\nu) &= \sup_{f^\pi\in\mathcal{F}_\pi}\left| \int f^\pi d\mu - \int f^\pi d\nu \right| \\
    & = \max_{\ell}\sup_{f\in\mathcal{F}} \left| \int \pi_\ell \circ f d\mu - \int \pi_\ell \circ f d\nu \right| \\
    & = \max_{\ell}\sup_{f\in\mathcal{F}} \left| \int f_\ell d\mu - \int f_\ell d\nu \right| \\
    & = \max_{\ell}\sup_{f_\ell\in\mathcal{F}_\ell} \left| \int f_\ell d\mu - \int f_\ell d\nu \right| \\
    & = \max_\ell \gamma_{\mathcal{F}_\ell}(\mu,\nu) = 0.
\end{align*}
Conversely, assuming $\gamma_{\mathcal{F}_\pi}(\mu,\nu) = 0$ yields
\begin{align*}
    \gamma_{\mathcal{F}_i}(\mu,\nu) & = \sup_{f_\ell\in\mathcal{F}_\ell} \left| \int f_\ell d\mu - \int f_\ell d\nu \right| \\
    & = \sup_{f\in\mathcal{F}} \left| \int \pi_\ell \circ f d\mu - \int \pi_\ell \circ f d\nu \right| \\
    & \leq \max_{\ell}\sup_{f\in\mathcal{F}} \left| \int \pi_\ell \circ f d\mu - \int \pi_\ell \circ f d\nu \right| \\
    & = \gamma_{\mathcal{F}_\pi}(\mu,\nu) = 0.
\end{align*}
This proves the first two equivalences. The third one is a byproduct of the proof.
\end{proof}

\subsection*{Proof of Corollary \ref{corollary:lpnorm}}
\begin{proof}
To avoid cumbersome notation, let $v = \frac{1}{n_0}\sum_{j=1}^{n_0} f^*(X^0_j) - \frac{1}{n_1}\sum_{j=1}^{n_1} f^*(X^1_j)$ and note $v_\ell = \frac{1}{n_0}\sum_{j=1}^{n_0} f_\ell^*(X^0_j) - \frac{1}{n_1}\sum_{j=1}^{n_1} f_\ell^*(X^1_j)$, then:
\begin{align*}
    Pr\left( \left\Vert v \right\Vert_{\ell_p}>\delta\right) & = Pr\left( \left\Vert v \right\Vert_{\ell_p}^p>\delta^p\right) \\
    & =
    Pr\left(\sum_\ell \left|v_\ell\right|^p > \delta^p \right) \\
    &\leq Pr\left(\sum_\ell \gamma_{\mathcal{F}_\ell}(\mathbb{Q}_{n_0}^0,\mathbb{Q}_{n_1}^1)^p > \delta^p\right) \\
    & \leq \sum_\ell Pr\left( \gamma_{\mathcal{F}_\ell}(\mathbb{Q}_{n_0}^0,\mathbb{Q}_{n_1}^1)^p > \delta^p/d\right) \\
    & = \sum_\ell Pr\left( \gamma_{\mathcal{F}_\ell}(\mathbb{Q}_{n_0}^0,\mathbb{Q}_{n_1}^1) > \delta/d^{1/p}\right) \\
    & \leq \sum_\ell B(\delta/d^{1/p},D^*,C^*) = dB(\delta/d^{1/p},D^*,C^*),
\end{align*}
where the second and third inequalities follow from a slight variation of Lemma~\ref{lemma:gammasum} and application of Lemma~\ref{lemma:fequiv}. For the $\ell_{\infty}$ case we have:
\begin{align*}
    Pr\left( \left\Vert v \right\Vert_{\ell_\infty}>\delta\right) & \leq Pr\left(\max_\ell \left|\gamma_\ell\right| > \delta \right)\\
    & \leq \sum_\ell B(\delta,D^*,C^*),
\end{align*}
concluding the proof.
\end{proof}

\subsection*{Balance for coarsening functions}
We will show the coarsened exact matching procedure belongs to a class of functions with tractable Vapnik-Chervonenkis dimension. Consider the set $\mathcal{S}$ of partitions with a fixed number of elements $R$. For a given partition $S\in\mathcal{S}$, such that $S = \{s_1, \dots,s_R\}$ define $f_S^{k\alpha}$ to be:
$$
f_S^{k\alpha}(x) = \sum_{i=1}^R k_i\alpha_i\chi_{s_i}(x),
$$
where $k_i\leq k$ for $k$ a constant, $\chi_{s_i}$ is the indicator function of $s_i$, and $\alpha:=(\alpha_1,\dots,\alpha_R)$ is a binary vector, this is, $\alpha_i\in\{0,1\}$ for each $i$. In words, if $x$ is found in $s_i$, $f$ will return a scaled version of $x$ if $\alpha_i$ is $1$ and zero otherwise.

Now let $\mathcal{F}:=\{f_S^{k\alpha}\}_{S\in\mathcal{S},\alpha\in A, k\leq \kappa}$, where $A$ is the set of all binary vectors of size $R$ and $\kappa\in \mathbb{R}$.
Hence, the coarsened exact matching procedure belongs to this class of functions, since in that case $\alpha_i$ indicates if there are at least two members of different groups in stratum $s_i$. For any sample point $x$, the weights are usually chosen in the following manner: If $x$ is a treated unit, $w_i^1 = 1$, otherwise, $w_i^0 = (m_1^{s}/m_1) / (m_0^s/m_0)$, where $s$ is the stratum $x$ belongs to. Letting $k_i = w_i^{\ell}n_\ell/m_\ell$ appropriately weighs matched samples.
We just need to add the mild assumption that the ratio of sample to matched size per stratum $s$ does not grow faster than $\sqrt{\kappa}$, that is, $n_\ell/m_\ell^s\leq \sqrt{\kappa}$ for all $s\in S$, because in that case $w_i^0\leq m_0/m_0^s\leq n_0/m_0^s\leq \sqrt{\kappa}$ and $n_\ell/m_\ell\leq \sqrt{\kappa} m_\ell^s/m_\ell\leq \sqrt{\kappa}$, so $k_i\leq \kappa$. Finally, notice that any similar function with a smaller partition size can be expressed by a function in $\mathcal{F}$, so we can consider variable partition size as long as it does not exceed a reasonable bound $R$.

For any set of points of size $R$ there is a partition $S$ containing one point in a different element, and therefore an $\alpha$ that can assign each point arbitrarily to either $0$ or $1$. So $\mathcal{F}$ shatters such set. However, if we add an extra point, and since the number of partitions is constrained, it would have to share partition element with a previous point, and so assignment under $f_s^{k\alpha}$. So the Vapnik-Chervonenkis dimension of $\mathcal{F}$ is $R$. Finally, let $g(\mathcal{Z}_\ell) = \mathbb{Q}_{n_\ell}^\ell$, where $\mathbb{Q}_{n_\ell}^\ell$ is the empirical distribution of the sample $\mathcal{Z}_\ell$ for group $\ell$. Let $k^*$ be chosen as above and let $(S^*,\alpha^*)$ be the particular partition and binary vector used for coarsened exact matching. Then, for the $\ell^{th}$ component we get:
\begin{align*}
    \left|\frac{1}{m_0}\sum_{i\in M_0}w_i^0 Z_{i,\ell}^0 - \frac{1}{m_1}\sum_{j\in M_1}w_j^1 Z_{j,\ell}^1 \right| &= \left|\frac{1}{n_0}\sum_{i=1}^{n_0}f_{S^*,\ell}^{k^*\alpha^*}(Z_i^0) - \frac{1}{n_1}\sum_{j=1}^{n_1}f_{S^*,\ell}^{k^*\alpha^*}(Z_j^1)\right| \\
    & \leq \sup_{f_{\ell}\in\mathcal{F^*}}\left|\frac{1}{n_0}\sum_{i=1}^{n_0}f_{\ell}(Z_i^0) - \frac{1}{n_1}\sum_{j=1}^{n_1}f_{\ell}(Z_j^1)\right| \\
    & = \gamma_{\mathcal{F}^*}(\mathbb{Q}^0_{n_0},\mathbb{Q}^1_{n_1}) = \gamma_{\mathcal{F}^*}(g(\mathcal{Z}_0),g(\mathcal{Z}_1)).
\end{align*}

Thus, the discrepancy among the matched samples per dimension is bounded by the $\gamma_{\mathcal{F}^*}$ distance of the unmatched samples.
Finally, the function $h(x):=\kappa x$ is an envelope function of $\mathcal{F}$ and has norm $\| h \|_{L_2(\mu)} < \infty$ as long as we assume compact domain, which is OK to do for most coarsened exact matching cases. Then, by Theorem 2.6.7 of \citet{van1996weak}:
\begin{align*}
    \sup_{\mu}  N(\epsilon,{\mathcal{F}},L_2(\mu)) \leq \left ( \frac{K}{\epsilon} \right)^{C^*},
\end{align*}
for some constant $K$ and where $C^*=2(R-1)$.

This leads us to our final result: Assume ideal balance on the population probabilities holds for $\gamma_{\mathcal{F}_\pi}$, then, for the $\ell^{th}$ component we have:
\begin{equation*}
    Pr\left(\left|\frac{1}{m_0}\sum_{i\in M_0}w_i^0 Z_{i,\ell}^0 - \frac{1}{m_1}\sum_{j\in M_1}w_j^1 Z_{j,\ell}^1 \right| > \delta\right) \leq B(\delta,D,C^*).
\end{equation*}
If we are interested in the $\ell_p$ norm of the full vector instead, then, by Corollary~\ref{corollary:lpnorm}:
\begin{equation*}
    Pr \left\{\left\Vert\frac{1}{m_0}\sum_{i\in M_0}w_i^0Z_i^0 - \frac{1}{m_1}\sum_{j\in M_1}w_j^1Z_j^1 \right\Vert_{\ell_p} > \delta\right\}     \leq dB(\delta/d^{1/p},D,C^*),
\end{equation*}
for finite $p\geq 1$. While
\begin{equation*}
    Pr \left\{\left\Vert\frac{1}{m_0}\sum_{i\in M_0}w_i^0Z_i^0 - \frac{1}{m_1}\sum_{j\in M_1}w_j^1Z_j^1 \right\Vert_{\ell_\infty} > \delta\right\}     \leq dB(\delta,D,C^*).
\end{equation*}

\subsection*{Balance using propensity scores}
Recall $e(Z) = P(T=1 \mid Z)$, and that we are assuming $Z \mid T=\ell \sim N(\mu_{\ell},\Sigma)$. Let $p_{\ell}$ be the probability density function of $N(\mu_{\ell},\Sigma)$, that is, the gaussian density, then by the density version of Bayes' Theorem we have
\[
p(T=1\mid Z=z) = \frac{p_1P(T=1)}{p_1P(T=1) + p_0P(T=0)}.
\]
Therefore, we can express the \textit{logit} of $e(Z)$ as
\[
logit(e(Z)) = \log\left(\frac{e(Z)}{1 - e(Z)}\right) = \log\left(\frac{p_1P(T=1)}{p_0P(T=0)}\right).
\]
Now define $L_k := logit(e(Z_k))$, then the matching procedure is based on the difference $\left| L_i - L_j \right|$. Given the above computation and after a few straightforward steps we get
\begin{align*}
    \left|
    L_i - L_j
    \right|
    & = \left|
    (\mu_1 - \mu_0)^T\Sigma^{-1}(Z_i - Z_j)
    \right| \\
    & = \left|
    f^*(Z_i) - f^*(Z_j)
    \right|,
\end{align*}
where $f^*(x) = w^T x$ for $w\in \mathbb{R}^{p}$. Notice the vector $w$ is the same as the one used for linear discriminant analysis so, adding an offset parameter, it will be useful to think of $f^*$ as a hyperplane.

Let $M_0^j$ be the control units assigned to treatment unit $j$. We make the assumption that there is a fixed number of assigned controls to each treatment, and so $m_0 = |M_0^j| m_1$. Then
\begin{align*}
    \Delta & := \left| \frac{1}{m_1}\sum_{j\in M_1} logit(e_j) 
    - \frac{1}{m_0}\sum_{i\in M_0} logit(e_i) \right| \\
    & = \left| \frac{1}{m_1}\sum_{j\in M_1} L_j 
    - \sum_{j\in M_1}\frac{1}{m_0}\sum_{i\in M_0^j} L_i \right| \\
    & = \left| \sum_{j\in M_1} \left(
    \frac{1}{m_1}L_j - \frac{1}{m_0}\sum_{i\in M_0^j} L_i
    \right)
    \right| \\
    & = \left| \sum_{j\in M_1} \left(
    \frac{1}{m_1}\sum_{i\in M_0^j}\frac{L_j}{|M_0^j|} - \frac{1}{m_0}\sum_{i\in M_0^j} L_i
    \right)
    \right| \\
    & = \left|
    \sum_{j\in M_1} \sum_{i\in M_0^j}
    \left(
    \frac{L_j}{m_1|M_0^j|} - \frac{L_i}{m_0}
    \right)
    \right| \\
    & = \left|
    \sum_{j\in M_1} \sum_{i\in M_0^j} \frac{1}{m_0}
    \left( L_j - L_i \right)
    \right| \\
    & = \left|
    \sum_{j\in M_1} \sum_{i\in M_0^j} \frac{1}{m_0}
    \left( f^*(Z_j) - f^*(Z_i) \right)
    \right| \\
    & = \left|
    \frac{1}{m_1}\sum_{j\in M_1}f^*(Z_j) - 
    \frac{1}{m_0}\sum_{i \in M_0}f^*(Z_i)
    \right|.
\end{align*}
That is, we can express the difference of means of $logit$s in terms of the difference of means of the discriminant functions. Let $p$ be the dimension of the covariates, and let $\mathcal{F}$ be the collection of $p$-dimensional hyperplanes, notice $f^*\in\mathcal{F}$. The Vapnik-Chervonenkis dimension of $\mathcal{F}$ is known to be $p+1$ \citep{mohri2018foundations}. We would like to bound $\Delta$ in terms of $\gamma$ but we first need some adjustments to $f^*$.

The matching procedure determines a set $\mathcal{Z}_M = \{Z_k \mid k\in M\}$ of matched samples, where $M =M_0\cup M_1$. By the Gaussian assumption the $Z$s are sampled from a Gaussian mixture so the probability of two sample points being the same is zero. Hence there is an $\epsilon>0$ such that for all $k\in M$, $\mathcal{Z}\cap B_\epsilon(Z_k) = \{Z_k\}$, that is, each $\epsilon$ ball centered around a matched sample does not contain any other sample point (here $\mathcal{Z}$ is the sample set). Let $S_\epsilon = \cup_k B_\epsilon(Z_k)$. Note $S_\epsilon$ is a measurable set. Let $\beta_{S_\epsilon}(x):= x \chi^{S_\epsilon}(x)$, this function maps points to zero if unmatched and to themselves if matched. Furthermore, let $\beta_\ell(x): = \frac{m_\ell}{n_\ell}\chi^{M_\ell}(x) + \chi^{M^C_\ell}$(x), for $\ell\in \{0,1\}$. Each $\beta_\ell$ scales elements in $M_\ell$ by the factor $\frac{m_\ell}{n_\ell}$ and leaves the rest untouched.

Notice $f_M^*:=f^*\circ\beta_1\circ\beta_0\circ\beta^{S_\epsilon}$ sends $Z_k$ to $\frac{m_\ell}{n_\ell}w^TZ_k$ if $k\in M_k$ and to $0$ otherwise. Then we can express $\Delta$ as
\begin{align*}
    \Delta & = \left|
    \frac{1}{m_1}\sum_{j\in M_1}f^*(Z_j) - 
    \frac{1}{m_0}\sum_{i \in M_0}f^*(Z_i)
    \right| \\
    & = \left|
    \frac{1}{n_1}\sum_{j=1}^{n_1}f_M^*(Z_j) - 
    \frac{1}{n_0}\sum_{i=1}^{n_0}f_M^*(Z_i)
    \right|.
\end{align*}

Now, consider the set $\mathcal{F}_M:=\{ f\circ\beta_1\circ\beta_0\circ\beta_{S} | f\in\mathcal{F}, S\in\Sigma\}$, where $\Sigma$ is the set of measurable sets according to the distribution of the $Z$s. The Vapnik-Chervonenkis dimension for $\mathcal{F}_M$ is the same as that of $\mathcal{F}$, that is, $p+1$. To see this we notice that the standard derivation for the hyperplane case involves shattering the standard basis $\mathcal{B}$ in $\mathbb{R}^p$. With probability one, no sample point will equal a standard basis vector, so there is an $\epsilon'>0$ for which we can create a set $s = \cup_{x\in\mathcal{B}}B_{\epsilon'}(x)$ such that $s\in\Sigma$ and no sample point is in $s$. Considering the functions $\{f_\nu\}$ in $\mathcal{F}$ used to shatter $\mathcal{B}$ and using $s$, we can use the functions $\{f_\nu\circ\beta_1\circ\beta_0\circ\beta_{s}\}$ in $\mathcal{F}_M$ to also shatter $\mathcal{B}$. So the Vapnik-Chervonenkis dimension is at least $p+1$. Since the functions $\beta_1$, $\beta_0$, and $\beta^S$ are either zero or a scaled identity, we don't get any complexity and the dimension is no larger than $p+1$, so it is indeed $p+1$. For the envelope function, we can choose $h(x) = <w_e,x>$. The norm of $w_e$ must be large enough to keep a $p+1$ Vapnik-Chervonenkis dimension. Since the vectors used to ensure such a dimension have norm $p+1$, the norm of $w_e$ must be at least $p+1$. So we can choose any large constant $C>p+1$. Since we are interested in vectors of the form $w=\Sigma^{-1}\Delta\mu$, we have $\|w\|\leq\|S^{-1}\|_{F}\|\Delta\mu\|_{2}$, so the user has to choose constants that bound each of these norms. Also, we must assume the covariates themselves are bounded, this ensures a finite norm for $h$.

Finally, we have
\begin{align*}
    \Delta & = \left|
    \frac{1}{n_1}\sum_{j=1}^{n_1}f_M^*(Z_j) - 
    \frac{1}{n_0}\sum_{i=1}^{n_0}f_M^*(Z_i)
    \right| \\
    & \leq \sup_{f\in \mathcal{F}_M} \left|
    \frac{1}{n_1}\sum_{j=1}^{n_1}f(Z_j) - 
    \frac{1}{n_0}\sum_{i=0}^{n_0}f(Z_i)
    \right| \\
    & = \gamma_{\mathcal{F}_M}(\mathbb{Q}^0_{n_0},\mathbb{Q}^1_{n_1}).
\end{align*}
Assuming Ideal Balance on the population probabilities, and applying Theorem 2.6.7 of \citet{van1996weak} in conjunction with Theorem~\ref{thm:concentration_theorem}, yields
\begin{equation*}
    Pr \{ \Delta > \delta\} \leq B(\delta,D,2p).
\end{equation*}

\subsection*{Covering number bound for Reproducing Kernel Hilbert Spaces}
We refer the reader to \citet{Wahba1990spline, berlinet2011reproducing, steinwart2008support} for nice overviews on reproducing kernel Hilbert spaces. Roughly speaking, a mapping $k: \cal X \times \cal X \rightarrow \mathbb{R}$ is said to be the reproducing kernel associated to the reproducing kernel Hilbert space $\mathcal{H}$ if it satisfies the following properties: (a) $k(\cdot,x) \in {\cal H}$ for any $x \in {\cal X}$; (b) $f(x) = \langle f,k(\cdot,x)\rangle_{\cal H}$ for all $f \in {\cal H}$ and $x \in {\cal X}$. Property (b) is commonly referred to as the reproducing property.

To apply Theorem~\ref{thm:concentration_theorem} to the reproducing kernel case, we will need to directly bound the covering number based on arguments different from Vapnik-Chervonenkis theory. Define the space
\[
{\cal H}^{m}_q(\mathbb{R}^p) = \{ f \in L_q(\mathbb{R}^p) \,\, | \,\, D^jf \in L_q(\mathbb{R}^p) \,\,\, \forall j\in\{1,\ldots,m\}; \,\,\, \|f\|_{q} < \infty \},
\]
where
\[
\|f\|_q = \sum_{0 \leq |\alpha| \leq s} \|D^{\alpha} f\|_{L_q}
\]
and $D^{\alpha}$ denotes partial derivatives in the sense of distributions. Then as a consequence of Theorem 1 of \citet{nickl2007bracketing}, if $m - q/p > 0$, then

\[ N(\epsilon,{\cal H},\|\cdot\|_q) \leq b_1 \epsilon^{-q},\]
while if $m - q/p < 0$,
\[ N(\epsilon,{\cal H},\|\cdot\|_q) \leq b_2 \epsilon^{-p/m},\]
Based on this result, Theorem~\ref{thm:concentration_theorem} can then be applied to prove a convergence rate under ideal balance. Note that this does not cover the Gaussian kernel case, because the Gaussian kernel is infinitely differentiable, so the space ${\cal H}^m_q(\mathbb{R}^p)$ does not apply. For the reader interested in the Gaussian case, we refer them to the recent paper by \cite{steinwart2020closer}.


\bibliographystyle{chicago}
\bibliography{biblio-IB.bib}

\end{document}